\documentclass[11pt]{article}
\usepackage{amsfonts}
\usepackage{amssymb}
\usepackage{amscd}
\usepackage{amsmath}
\usepackage{epsfig}
\usepackage{graphicx, subfigure}
\usepackage{latexsym}
\usepackage{pst-3dplot}
\usepackage[symbol*]{footmisc}
\newlength{\dinwidth}
\newlength{\dinmargin}
\usepackage[all,cmtip]{xy}

\newtheorem{definition}{Definition}
\newtheorem{theorem}{Theorem}

\newtheorem{proposition}{Proposition}
\newtheorem{corollary}{Corollary}
\newtheorem{remark}{Remark}
\newtheorem{lemma}{Lemma}
\newtheorem{example}{Example}

\def \P{\mathbb P}

\begin{document}

\renewcommand*{\thefootnote}{\fnsymbol{footnote}}

\title{A new relationship between block designs}
\author{A. Shramchenko, V. Shramchenko$^*$}

\date{}

\maketitle

\footnotetext[1]{Department of mathematics, University of
Sherbrooke, 2500, boul. de l'Universit\'e,  J1K 2R1 Sherbrooke, Quebec, Canada. E-mail: {\tt Vasilisa.Shramchenko@Usherbrooke.ca}}

 \begin{abstract}
We propose a procedure of constructing new block designs starting from a given one by looking at the intersections of its blocks with various sets and grouping those sets according to the structure of the intersections. 
 We introduce a symmetric relationship of friendship between block designs built on a set $V$ and consider families of block designs where all designs are friends of each other, the so-called friendly families. We show that a friendly family admits a partial ordering. Furthermore, we exhibit a map from the power set of $V$, partially ordered by inclusion, to a friendly family of a particular type which preserves the partial order. 
 
  \end{abstract}
\maketitle

\section{Introduction}

We consider balanced incomplete block designs (BIBDs), see \cite{encyclopedia, handbook, Rib}.
With every block design one can associate some other block designs obtained using the original one. For example, given a block design $D$ with blocks' length $k$, one can form another block design by taking all subsets of length $k'<k$ of all the blocks of $D$. We suggest a new way of looking at relationships between block designs. 

We introduce a notion of {\it friendship} between block designs based on the structure of intersections between their blocks. More precisely, suppose we choose a block of one of the two BIBDs and look at the intersections of this block with all the blocks of the other BIBD. Suppose also that the number of blocks of the second BIBD that give intersections of length $n$ only depends on $n$ and does not depend on our choice of a block of the first BIBD. Now, if the same is true after we interchange the roles of the two BIBDs, then they are called {\it friends}. This is quite a strong condition and we expect families of block designs that are friends of each other to possess interesting properties. We give examples of friends and study the question of when a block design is friends with itself. 


We call a family of block designs built on a set $V$ which are pairwise friends a {\it friendly} family. On such a family we 
introduce a partial order. Furthermore, given that the power set of $V$ is also partially ordered by inclusion, we show that there exists a map from the power set of $V$ to a friendly family of a particular type which preserves the ordering.

We also suggest to study families of block designs constructed using a common ``parent" block design. In Section \ref{sect_construction}, making use of a given block design built on a set $V$, we define an equivalence relation on the power set of the set $V$. 
We conjecture that if the given block design is a Desargues projective plane, then the constructed equivalence classes in the power set are block designs, which are pairwise friends, that is form a friendly family.
Furthermore, we suggest to use the procedure of intersecting the blocks of a BIBD by various sets as a way to generate new block designs. We describe a few cases when such a construction yields friends of the original block design.

\section{Definitions and notation}

A balanced incomplete block design is defined as follows, see \cite{encyclopedia, handbook, Rib}.

\begin{definition}
Let $V$ be a finite set of cardinality $v=|V|$ and $b,k,r,\lambda$ four positive integers. A block design with parameters $(v,b,r,k,\lambda)$ built on $V$ is a list of $b$ blocks, each of which is a $k$-element subset of $V$, such that every element of $V$ is contained in exactly $r$ blocks and every pair of elements of $V$ is contained in exactly $\lambda$ blocks.
\end{definition}

 We assume that all block designs are simple, i.e. have no repeated blocks. The notation we use is in accordance with \cite{Rib}.

\begin{itemize}


 \item ${\mathcal P}(V)$ the power set of $V$, that is the set of all subsets of $V$ including the empty set and $V$ itself, ordered by inclusion: for $X,Y \subset V$ we say $X<Y$ if and only if $X\subset Y$.

 \item  $D_i, i \in I$, a set of block designs, each of which is built on the set $V$ and has parameters 
 $(v,b_i,r_i, k_i,\lambda_i)$. 

 \item $D^s_i, s=1, \dots,  b_i$, stand for the blocks of the block design $D_i$.


  
\item $\mathcal D_k$ is a full design of block size $k$ built on the set $V$, that is a block design with parameters $(v, {v \choose k},r,k,\lambda)$ whose set of blocks is the full set of combinations of $k$ out of $v$ elements. We assume that $D_0$ is the empty set and $\mathcal D_v=V.$   
  
 \end{itemize}

\begin{remark}
\label{rmk_relations}
The parameters of a block design are not independent: $\!$one has $bk\!=\!vr\!$ and $r(k-\!1)\!=\!\lambda(v-1\!)$.
\end{remark}

\section{Friends of block designs}

We are going to consider intersections of the blocks of a design by a given set and focus on the number of elements in these intersections.  Let us consider a block design $D$ built on a set $V$ and some subset $M$ of $V$. For $j\geq 0$, let  $z_j\in\mathbb N$ stand for a number of blocks of $D$ whose intersections with $M$ contain exactly $j$ elements. More precisely, we define
\begin{equation}\label{mama} 
{ z_j=\displaystyle\sum\nolimits_{1\leq s\leq b} \delta( |M \cap B_s|,j )},
\end{equation}
where for two sets $M$ and $L$
\begin{equation*}
 \delta( |M \cap L|,j )=
 \begin{cases}
  1& \text{if $|M \cap L|=j$,}\\
  0& \text{otherwise.}
 \end{cases} 
\end{equation*}

\begin{lemma} \label{lemma_formulas3}
Let ${D}$ be a  block design with parameters ${(v,b,r,k,\lambda)}$ built on the set $V$. Let $M$ be a subset of $V$ with $|M|=m$. Then the following  formulas hold:
\label{lemma_formulas3}
\begin{eqnarray}
\label{eq1}
&&\sum_{j=0}^k z_j=b,\\
\label{eq2}
&&\sum_{j=0}^k z_j{j}=rm,\\
\label{jsquare}
&&\sum_{j=0}^k z_j{j^2}=m(\lambda{m}-\lambda+r).
\end{eqnarray} 
\end{lemma}
{\it Proof.}
The first two equations follow directly from the definition of a block design. In the same way one obtains
 \begin{equation}
\sum_{j=0}^{ k} z_j {j \choose 2} = {m \choose 2}  \lambda .
  \end{equation} 
  Here and below we assume that ${m \choose n}=0$ if $n>m.$
This and the first two equations of the lemma imply (\ref{jsquare}).
$\Box$

Due to the first equation of Lemma \ref{lemma_formulas3}, the sequence ${\varphi= (z_0, z_1, . . . , z_k)}$ is a partition of the number of blocks  $b$. Let us denote this partition by $\varphi(D,M)=(z_0, z_1, . . . , z_k).$ Note that we allow $z_i=0$ for some $i \leq k$, and do not require $z_i \geq z_j$, for   $i \geq j$.

 \begin{definition}
Two block designs $D_1$ and $D_2$ built on the same set $V$ are called friends if the partitions $\varphi(D_1,D_2^i)$ and $\varphi(D_2,D_1^j)$ do not depend on $i$ and $j$, respectively. 
\end{definition}

Note that this relationship between designs is symmetric by definition. 

For designs which are friends, we simplify the notation for partition: we write $\varphi(D_1,D_2)$ instead of $\varphi(D_1,D^i_2)$ and $\varphi(D_2,D_1)$ instead of $\varphi(D_2,D^i_1)$. 

\begin{example}
\label{P3}
{\rm
Let  $\P_3$ be the finite projective plane with parameters $(v=b=7, r=3, k=3, \lambda=1)$ built on the set $V=\{1,2,\dots,7\}$, the Fano plane. Denote its blocks by $B_i$, that is $\mathbb P_3=\{B_1,\dots, B_7\}$.
More precisely, we have $B_1=(2,3,5), \; B_2=(3,4,6),\; B_3=(4,5,7), \;B_4=(1,5,6), \;B_5=(2,6,7),\; B_6=(1,3,7), \;B_7=(1,2,4).$
Let $\mathcal D_5$ be the full design of block size five built on the same set $V$.
Then $\P_3$ and $\mathcal D_5$ are friends. The corresponding partitions are $\varphi(\mathcal D_5, \P_3)=(0,3,12,6)$ and $\varphi(\P_3, \mathcal D_5)=(0,1,4,2).$

}
\end{example}

For more examples of friends see Section \ref{sect_friendly}.

 \begin{proposition}
  Let block designs $D_1$ and $D_2$ be friends, then:
   \begin{equation}\label{unic}
    {\varphi(D_1,D_2)b_2=\varphi(D_2,D_1)b_1}.
   \end{equation}
 \end{proposition}
{\it {Proof.}}
Consider a matrix $\mathcal M$ with entries given by $\mathcal M_{ij}=|D_1^i \cap D_2^j| $ for  $i=1,\dots,b_1$ {and }    $j=1,\dots,b_2$. Because $D_1$ and $D_2$ are friends, the rows (columns)  of $\mathcal M$ are equal up to permutations. Any integer $n$ appears equal number of times in every row (column). Multiplying this number of times by $b_1$ (respectively $b_2$), we get
the total number of times the integer $n$ appears in the matrix $\mathcal M$. On the other hand, by doing so we get the corresponding component of the partition in the left (respectively, right) hand side of   \eqref{unic}.
$\Box$

\begin{proposition}
\label{prop_complement}
Let block designs $D$ and $D_1$ be friends with the corresponding partition $\varphi(D, D_1) = (z_0, z_1, \dots, z_{k_D})$; here ${k_D}$ is the block size of $D$. Let design $D_2$ be the complement of $D_1$, i.e. $D^s_2=V\setminus D^s_1$ for all $s=1,\dots, b_1$. Then $D$ and $D_2$ are also friends and the  partition $\varphi(D, D_2)=(\omega_0,\dots, \omega_{k_D})$ satisfies: 
\begin{equation*}
\omega_i=z_{k_D-i}, \qquad i=0,\dots, k_D.
\end{equation*}
\end{proposition}

{\it Proof.} Let us assume that $k_D$ is smaller than block sizes of designs $D_1$ and $D_2$; the other situations are considered similarly. Suppose a block $D^s_1$ has $i$ elements in common with a block $D^t$ of design $D$. Then the remaining $k_D-i$ elements of $D^t$ belong to the complement of $D^s_1$ that is to $D^s_2$. Thus $\omega_i=z_{k_D-i}$.
$\Box$

\vskip 0.5cm
A natural question is when a block design is friends with itself. 

\begin{example}{\rm
A full block design $\mathcal D_k$ is friends with itself. Indeed, in this case  ${\varphi(\mathcal D_k,\mathcal D_k^j)}=(z_0, z_1, . . . , z_k)$ is independent of the choice of a block of  $\mathcal D_k$ since for all such partitions we have
  \begin{equation}
  {z_i=
  \begin{pmatrix}
k\\i
\end{pmatrix}
  \begin{pmatrix}
v-k\\k-i
\end{pmatrix}}.
  \end{equation}
}
\end{example}

Here are some classes of block designs that are naturally friends of themselves. 

\begin{theorem}
\label{prop_conj}
Let $D$ be a block design with parameters $(v,b,r,k,\lambda)$. 
	\begin{enumerate}
	\item If $\lambda=1$, then $D$ is friends with itself. 
	\item If $k=3$, then $D$ is friends with itself. 
	\item If $D$ is symmetric, that is $b=v$, then $D$ is friends with itself. 
	\end{enumerate}
\end{theorem}
{\it Proof.} 
\begin{enumerate}
\item Let $B$ be an arbitrary block of $D$. As before, we denote $\varphi(D, B)=(z_0, z_1, . . . , z_k)$ a partition of the number $b$ of blocks obtained by intersecting all blocks of $D$ with the chosen block $B$. 

Given that every pair of elements appears in only one block of $D$ ($\lambda=1$), we know that no pair of elements from the block $B$ will appear in some other block of $D$. Thus there will be no blocks whose intersection with $B$ would give a set of two or more elements with the only exception of the set of  $k$ elements produced by the intersection of $B$ with itself. Thus $z_i=0$ for $2\leq i\leq k-1$ and $z_k=1$.  

The remaining elements $z_0$ and $z_1$ can be found from equations of Lemma \ref{lemma_formulas3}. 
Thus the partition $\varphi(D, B)$ is independent of the choice of a block $B$.

\item Let us first consider a block design $D$ with an arbitrary $k$. For some arbitrary block $B$ of $D$ denote $\varphi(D, B)=(z_0, z_1, . . . , z_k)$. Lemma \ref{lemma_formulas3} with $m=k$ gives three relations for $\{z_{j}\}_{j=0}^k$.
%
Our condition that no blocks are repeated implies $z_k=1$. Thus we have four linear equations for $k+1$ partition elements $z_0,\dots, z_k.$ Therefore, when $k=3$ we can find the partition starting from parameters of the block design. This means in particular, that the partition $\varphi(D, B)$ does not depend on the choice of a block $B$.

\item This is a simple corollary of the well known fact, see \cite{handbook} II.6, that in a symmetric design every two distinct blocks have $\lambda$ points in common. 
 $\Box$
\end{enumerate}

\begin{corollary} Any finite projective plane is friends with itself. 
\end{corollary}
{\it Proof.} Recall that a finite projective plane is a block design with $b=v$ and $\lambda=1$. $\Box$

\vskip 0.2cm

\begin{example}{\rm
Here is a design that is friends with itself. The parameters are $(v=9, b=12, r=8, k=6, \lambda=5)$ so this example is not covered by Theorem \ref{prop_conj}.
This design is taken from \cite{CC}, page 474. 
\begin{eqnarray*}
\left.\begin{array}{cccccc}1 & 2 & 4 & 5 & 7 & 8 \\2 & 3 & 5 & 6 & 8 & 9 \\1 & 3 & 4 & 6 &7 & 9\end{array}\right. 
\qquad\qquad
\left.\begin{array}{cccccc} 1& 3 & 5 & 6 & 7 & 8 \\1 & 2 & 4 & 6 & 8 & 9 \\2 & 3 & 4 & 5 & 7 & 9\end{array}\right.
\\
\\
\left.\begin{array}{cccccc}1 & 2 & 5 & 6 & 7 & 9 \\1 & 3 & 4 & 5 & 8 & 9 \\2 & 3 & 4 & 6 & 7 & 8\end{array}\right.
\qquad\qquad
\left.\begin{array}{cccccc}4 & 5 & 6 & 7 & 8 & 9 \\1 & 2 & 3 & 4 & 5 & 6 \\1 & 2 & 3 & 7 &8 & 9\end{array}\right.
\end{eqnarray*}
The corresponding partition of $12$ is $(0,0,0,2,9,0,1)$.
}
\end{example}

Although it is natural to suggest that every block design is friends with itself, this is not true. For block designs which are not friends with themselves, see Example \ref{example_sts}.

Transitivity does not hold for the relationship of friendship either. To see this we first prove the following simple lemma. 

\begin{lemma}
Let $D$ be a block design $(v,b,r,k,\lambda)$ built on a set $V$ such that $k<v-1.$ Then $D$ is friends with the full design $\mathcal{D}_{v-1}$.
\end{lemma}
{\it Proof.} It is straightforward to compute the partitions $\varphi(D, \mathcal{D}_{v-1}^i)=(z_0,z_1,\dots, z_k)$ and $\varphi(\mathcal{D}_{v-1}, D^i)=(w_0,w_1,\dots, w_k)$. One obtains $z_{k-1}=r$, $z_k=b-k$, $w_{k-1}=k$ and $w_k=v-k$; all other entries vanish. 
$\Box$

\begin{corollary}
The relationship of friendship is not transitive.
\end{corollary}
{\it Proof.} Consider two block designs built on the same set $V$ which are not friends of each other and such that their block sizes are smaller than $v-1$. By the lemma they are both friends with $\mathcal D_{v-1}$.
$\Box$

\section{Friendly families of block designs}
\label{sect_friendly}
In this section we consider a family of  designs which are pairwise friends of each other. Let us call such a family (or set) {\it friendly}. For example, the set of full  designs $\{\mathcal D_j\}_{j=0}^{v}$ is a friendly family.
Such a set admits a partial order as follows.

 \begin{definition}
Let two block designs $D$ and $\overline{D}$  be friends with   $k<\overline{k}$, where  $k$ and $\overline{k}$ are the  sizes of blocks of $D$ and $\overline{D}$, respectively, and $\varphi(D,\overline{D})=(z_0, z_1, . . . , z_k)$. We say that  $D<\overline{D}$ if $ z_k>0$.
\end{definition}
 
We thus have two partially ordered sets: a friendly set of block designs built on a set $V$ and the power set ${\mathcal P}(V)$ (ordered by inclusion). It turns out that there exists a map between these two sets that preserves the ordering.

\begin{proposition}
\label{prop_order}
Let $\mathbf{D}=\{D_i\}_{i\in I}$ be a friendly family of designs built on $V$.
 Suppose that no two designs of $\mathbf D$ share a block and that the set of all blocks of all designs in $\mathbf D$ gives the power set  ${\mathcal P}(V)$. Then there exists a map
 $\alpha: {\mathcal P}(V)  \rightarrow \mathbf{D}$  preserving the partial order.
\end{proposition}
\it {Proof.} {\rm The map
 $\alpha: {\mathcal P}(V)  \rightarrow \mathbf{D}$ is defined as follows. It sends a subset $U$ of $V$ to the  design in $\mathbf D$ which contains this subset $U$ as a block. By the assumptions of the proposition, there exists a unique block design for which $U$ is a block.  
Note that we consider the empty set and the full set $V$ 
as degenerate designs included in $\mathbf D$.

Now, suppose $X,Y \in {\mathcal P}(V)$,   $X \subset Y$,  $X \neq  Y$, that is $X< Y$. We want to show that $\alpha(X)< \alpha(Y)$. Since $|X| \neq  |Y|$, we have that $X$ and $Y$ belong to two distinct block designs from $\mathbf{D}$, say, $X \in D$, $Y \in \overline{D}$ with $k<\overline{k}$, where   $k$ and  $\overline{k}$ are block sizes of $D$ and $\overline{D}$ respectively. Since  $D$ and $ \overline{D}$ are friends and  $X< Y$, then for $\varphi(D,\overline{D})=(z_0, z_1, . . . , z_k)$ we have $ z_k>0$, and therefore $\alpha(X)< \alpha(Y)$.
}$\Box$

\vskip 0.3cm
%
{\rm


In the following example, we construct a friendly family of designs that satisfies conditions of Proposition \ref{prop_order} and that is  different from the set of full  designs $\{\mathcal D_j\}_{j=0}^{v}$.

\begin{example}
\label{k3}
{\rm 
All designs are built on the set $V=\{1,2,\dots,7\}.$

Let $\mathbb P_3=\{B_1,\dots, B_7\}$ be the finite projective plane from Example \ref{P3}.
Denote by $D'_3$ the  design with $b=28$ and $k=3$ whose blocks are given by all sets of three elements of $V$ which are not blocks of $\mathbb P_3.$

Denote by $D_4$ the  design with  $b=7$ and $k=4$ whose blocks are complements of the blocks of $\mathbb P_3$, that is $ D_4 = \{\bar{B}_1, \dots, \bar{B}_7\}$ where $\bar{B}_i=V\setminus B_i.$ Similarly, by $ D'_4$ we denote the  design with $b=28$ and $k=4$ whose blocks are the sets of four elements which are not blocks of $ D_4$.

To the family $\{\mathbb P_3, D'_3,D_4,D_4'\}$ we also add the empty set $\mathcal D_0$ and  $\mathcal D_1, \mathcal D_2, \mathcal D_5, \mathcal D_6, \mathcal D_7=V$  where $\mathcal D_k$ is the full design of block size $k$.

In this family all designs are pairwise friends. Moreover, we have $\mathcal D_1<\mathcal D_2$, $\mathcal D_2< \mathbb P_3$ and $\mathcal D_2<  D'_3$,  $ D'_3 <  D_4$ and $ D'_3 <  D'_4$, $\mathbb P_3<  D'_4$, $ D_4<\mathcal D_5$ and $ D'_4<\mathcal D_5<\mathcal D_6< \mathcal D_7$.
}
\end{example}

In this section we presented a new way of constructing a partially ordered set (a poset) given by a  full (satisfying conditions of Proposition \ref{prop_order}) friendly family of block designs. The condition of pairwise friendship seems to be a very strong one therefore we expect such families to possess interesting properties. There arise a few interesting questions, for example: to say how many different full friendly sets of block designs one can construct on a given set $V$, and to find analogues of the result of the papers \cite{BeckZ, posets} for such frienly families.



%

\section{Constructing friends of block designs}
\label{sect_construction}

In this section we use the procedure of intersecting blocks of a design by various sets to form block designs starting with a given one. 

Let $\mathbb D$ be a design built on a set $V$ whose blocks are of size $k$. 
Now, for every $n=1,2,\dots, v$ consider the set $\mathfrak{N}_n$ of all subsets of $V$ of size $n$. 

We now subdivide the set $\mathfrak{N}_n$ into classes of subsets $D_j^{(n)}$ as follows. 
Consider the map $\Phi$ from the set $\mathfrak{N}_n$ to the set of partitions of $v$ which for a set $S\in\mathfrak{N}_n$ gives $\Phi(S) = \varphi(\mathbb D, S).$ In this way, we obtain a number of partitions in the image: $\Phi(\mathfrak{N}_n)=\{\varphi^{(n)}_1, \dots, \varphi^{(n)}_{s_n}\}$. Denote by $D_j^{(n)}$ the class of sets from $\mathfrak{N}_n$ given by the inverse image of $\varphi^{(n)}_j$, that is $D_j^{(n)}=\Phi^{-1}(\varphi^{(n)}_j)$. 

The family of these classes forms a special subdivision of the power set $\mathcal P(V)$ into  non-intersecting equivalence classes.
In other words, we may say that two subsets of $V$ are equivalent if and only if they belong to the same class $D_j^{(n)}$ for some $n$ and $j$. 
Let us denote this subdivision of $\mathcal P(V)$ by $\mathfrak{M}$.

If the original  design is a projective plane on seven elements with blocks of size three, $\mathbb D=\P_3$, this subdivision coincides with the one described in Example \ref{k3}, that is  $\mathfrak{M}=\{ \mathcal D_0, \mathcal D_1, \mathcal D_2, \P_3, D'_3, D_4, D'_4,\mathcal D_5,\newline  \mathcal D_6,\mathcal D_7 \}$. For a projective plane with $k=4$ it is easy to verify that all the $D_j^{(n)}$ are block designs and that the family $\mathfrak{M}$ is  again friendly and satisfies conditions of Proposition \ref{prop_order}. 

We conjecture that if  the original  design $\mathbb D$ is a Desargues projective plane, the above construction yields a friendly family of designs.

The next theorem gives some cases when the described procedure results in block designs that are friends with the original one.

\begin{theorem}
\label{thm_friends}
Let $\mathbb D$ be a block design with parameters $(v,b,r,k=3,\lambda)$. Then 
\begin{itemize}
\item There are two classes of sets $D^{(3)}_1=\mathbb D$ and $D^{(3)}_2$ with the upper index $(3)$. Both of them are block designs and both are friends with $\mathbb D$.
\item If $\lambda=1$, there are two classes of sets $D^{(4)}_1$ and $D^{(4)}_2$ with the upper index $(4)$. In this case, both of them are block designs and both are friends with $\mathbb D$.
\end{itemize}
\end{theorem}
{\it Proof.}
\begin{itemize}
\item 
Define $D^{(3)}_2$ to be the set of all $3$-subsets of $V$ which are not blocks of $\mathbb D$, and define $D^{(3)}_1$ to coincide with $\mathbb D$. As is easy to see, all $3$-subsets of $V$ are covered by these two classes. Then $D^{(3)}_1=\mathbb D$ is trivially a block design and is friends with itself by Theorem \ref{prop_conj}. 

To prove that $D^{(3)}_2$ is a block design, we need to find its parameters $(v^{(3)}_2\; b^{(3)}_2, \;r^{(3)}_2,\; k^{(3)}_2,\; \lambda^{(3)}_2)$. Two of these parameters are known: $v^{(3)}_2=v$ and $ k^{(3)}_2=3$. The remaining parameters are obtained easily knowing that the blocks of $D^{(3)}_2$ and those of $\mathbb D$ exhaust all triples of elements of $V$, we have $b^{(3)}_2={v \choose 3}-b$, $r^{(3)}_2={v-1\choose 2}-r$ and  $\lambda^{(3)}_2=v-2-\lambda$.

Let us now prove that $\mathbb D$ and  $D^{(3)}_2$ are friends. 
 Let $B^i$ be a block of $D^{(3)}_2$ and consider the partition $\varphi(\mathbb D, B^i) = (z_0, z_1, z_2, z_3)$. By the definition of $D^{(3)}_2$ we have $z_3=0$. Three remaining $z_i$ are found from the three equations of  Lemma \ref{lemma_formulas3}. 
By switching the roles of  $\mathbb D$ and  $D^{(3)}_2$ in the above calculation, we obtain 
that $\varphi(D^{(3)}_2, \mathbb D^i)$  is also independent  of $i$. 

\item 

Define now $D^{(4)}_1=\{T\subset V,\; |T|=4, |\exists \; \mathbb D^i\subset T\}$ to be the set of all $4$-subsets of $V$ which contain at least one block of $\mathbb D$. And define $D^{(4)}_2=\{T\subset V,\; |T|=4, | \mathbb D^i\not \subset T \;\forall i\}$ to be the set of all $4$-subsets of $V$ which do not contain any block of $\mathbb D$. All $4$-subsets of $V$ are covered by these two classes.

Assuming $\lambda=1$,
let us first prove that $D^{(4)}_1$ and $D^{(4)}_2$ are block designs. We do this by computing their parameters $v^{(4)}_j,b^{(4)}_j,r^{(4)}_j,k^{(4)}_j,\lambda^{(4)}_j$ for $j=1,2$. We know that $v^{(4)}_j=v$ and $k^{(4)}_j=4$. 
To find $b^{(4)}_j$, let us note that the blocks of $D^{(4)}_1$ are obtained by taking a block of $\mathbb D$ and upending to it one element not already contained in the block. In this way, for every block of $\mathbb D$ we get $(v-3)$ new sets in $D^{(4)}_1$. Since $\lambda=1$,  sets obtained from different blocks of $\mathbb D$ will intersect by at most two elements. Thus they will all be distinct and we have $b^{(4)}_1=b(v-3).$ 
Now,  $r^{(4)}_1 = r(v-3)+(b-r)$. This is because for a given element for each of the $r$ blocks in $\mathbb D$ that contain it, we get $(v-3)$ blocks in $D^{(4)}_1$; similarly, for each of the $(b-r)$ blocks of $\mathbb D$ that do not contain the given element we can upend this element to obtain a block of $D^{(4)}_1$ that contains it. Reasoning analogously, we obtain $\lambda^{(4)}_1 = v-3+2(r-1).$
Parameters of  $D^{(4)}_2$ are determined knowing that the blocks of $D^{(4)}_1$ and $D^{(4)}_2$ exhaust all quadruples of elements of $V$.

Let us now prove that  $D^{(4)}_j$ and $\mathbb D$ are friends. We need to see that all four partitions 
$\varphi(D^{(4)}_j, \mathbb D^i)$ and $\varphi\left( \mathbb D, (D^{(4)}_j)^i\right)$ are independent of $i$. All of the partitions contain four parts $(z_0,z_1,z_2,z_3)$, thus by Lemma \ref{lemma_formulas3}, it is enough to determine one of the parts. As is easy to see,  if $\lambda=1$, the $z_3$ part of $\varphi(D^{(4)}_1, \mathbb D^i)$ is equal to $v-3$ and that of $\varphi(D^{(4)}_2, \mathbb D^i)$ vanishes. Similarly we find that the $z_3$ part of $\varphi\left( \mathbb D, (D^{(4)}_1)^i\right)$ is equal to one and that of $\varphi\left( \mathbb D, (D^{(4)}_2)^i\right)$ vanishes by definition. 
$\Box$
\end{itemize}

\begin{example}
\label{example_sts}
Consider the following two Steiner triple systems $STS(13)$, that is two block designs with parameters $(v=13,b=26,r=6,k=3,\lambda=1)$. We denote them by $S_1$ and $S_2$ and list all their blocks. The blocks are such that $S_1^i=S_2^i$ for $i=1,\dots, 22$, and the four remaining blocks are different in the two STS. Here is the list of the blocks. 
\begin{eqnarray*}
&&S_j^1=(1,2,3), \qquad S_j^2=(1,4,5), \qquad S_j^3=(1,6,7), \qquad S_j^4=(1,8,9), \qquad S_j^5=(1,10,11), \\
&&S_j^6=(1,12,13), \qquad S_j^7=(2,4,6), \qquad S_j^8=(2,5,7), \qquad S_j^9=(2,8,10), \qquad S_j^{10}=(2,9,12),\\
 && S_j^{11}=(2,11,13), \qquad S_j^{12}=(4,3,8), \qquad S_j^{13}=(4,7,9), \qquad S_j^{14}=(4,10,13), \\
 &&  S_j^{15}=(4,11,12),\qquad S_j^{16}=(7,3,11), \qquad S_j^{17}=(7,8,13), \qquad S_j^{18}=(7,10,12), \\
 && S_j^{19}=(8,5,11), \qquad S_j^{20}=(8,6,12),\qquad S_j^{21}=(6,9,11), \qquad S_j^{22}=(3,5,12)
\end{eqnarray*}
and the four remaining blocks in each STS are
\begin{eqnarray*}
S_1^{23}=(3,6,10), \qquad S_1^{24}=(3,9,13), \qquad S_1^{25}=(5,6,13), \qquad S_1^{26}=(5,9,10).
\end{eqnarray*}
\begin{eqnarray*}
S_2^{23}=(3,6,13), \qquad S_2^{24}=(3,9,10), \qquad S_2^{25}=(5,6,10), \qquad S_2^{26}=(5,9,13).
\end{eqnarray*}
On these two block designs we perform the procedure described in the beginning of this section and find the sets $D_j^{(n)}$ for $n=3,4,5,6.$ The corresponding sets for other values of $n$ can be obtained using Proposition \ref{prop_complement}.
\begin{itemize}
\item For both $S_1$ and $S_2$, the set $\Phi(\mathfrak N_3)$ contains two partitions, $\varphi_1^{(3)}=(10,15,0,1)$   and $\varphi_2^{(3)}=(11,12,3,0)$. The number of blocks in $D_1^{(3)}$ is $26$ and the number of blocks in $D_2^{(3)}$ is $260$.
\item For both $S_1$ and $S_2$, the set $\Phi(\mathfrak N_4)$ contains two partitions, $\varphi_1^{(4)}=(7,15,3,1)$ corresponding to $260$ blocks in $D_1^{(4)}$  and $\varphi_2^{(4)}=(8,12,6,0)$ corresponding to $455$ blocks in $D_2^{(4)}
$.
\item For both $S_1$ and $S_2$, the set $\Phi(\mathfrak N_5)$ contains three partitions: $\varphi_1^{(5)}=(5,13,7,1)$ corresponding to $780$ blocks in $D_1^{(5)}$, then $\varphi_2^{(5)}=(4,16,4,2)$ corresponding to $195$ blocks in $D_2^{(5)}$ and $\varphi_3^{(5)}=(6,10,10,0)$ corresponding to $312$ blocks in $D_3^{(5)}$.
\item For both $S_1$ and $S_2$, the set $\Phi(\mathfrak N_6)$  contains five partitions: $\varphi_1^{(6)}=(2,15,6,3)$, $\varphi_2^{(6)}=(1,18,3,4)$, $\varphi_3^{(6)}=(4,9,12,1)$, $\varphi_4^{(6)}=(3,12,9,2)$ and $\varphi_5^{(6)}=(5,6,15,0)$.
However the number of blocks in the corresponding sets $D_j^{(6)}$ differ for $S_1$ and $S_2$. Namely, the sets $D_1^{(6)}$ contain $208$ and $228$ blocks for $S_1$ and $S_2$, respectively. The numbers of blocks for $D_2^{(6)}$ are $13$ and $8$, for $D_3^{(6)}$ - $468$ and $488$, for $D_4^{(6)}$ - $988$ and $958$, for $D_5^{(6)}$ - $39$ and $34$, respectively.
\end{itemize}
 It turns out that for $n=3,4,5$ all the obtained sets $D_j^{(n)}$  are block designs, moreover  for $n=3$ and $n=4$, we obtain friendly families of block designs. For $n=5$ the obtained block designs are not friends of one another nor are they friends with themselves. They are, however, friends with the original Steiner triple system. For $n=6$ none of the sets $D_j^{(6)}$ is a block design. 
\end{example}

\vskip 0.5cm
{\bf Acknowledgments.} The second author gratefully acknowledges
support from the Natural Sciences and Engineering Research Council of Canada (discovery grant) and the University of Sherbrooke.


}

\end{document}